\documentclass[11pt,oneside,leqno]{amsart}%
\usepackage{amsmath}
\usepackage{amsfonts}
\usepackage{amssymb}
\usepackage{graphicx}
\usepackage{color}%
\providecommand{\U}[1]{\protect\rule{.1in}{.1in}}
\newtheorem{theorem}{Theorem}

\begin{document}
\title[Ricci solitons on low-dimensional generalized symmetric spaces]{Ricci solitons on low-dimensional generalized symmetric spaces}
\author{Giovanni Calvaruso and E. Rosado}
\address{Giovanni Calvaruso: Dipartimento di Matematica e Fisica \lq\lq E. De
Giorgi\rq\rq\\
Universit\`a del Salento\\
Prov. Lecce-Arnesano \\
73100 Lecce\\
Italy.\\
E. Rosado: {Department of Applied Mathematics, Escuela T\'{e}cnica Superior de
Arquitectura, Universidad Polit\'{e}cnica de Madrid}\\
Avda. Juan de Herrera 4\\
28040 Madrid \\
Spain.}
\email{{eugenia.rosado@upm.es}}
\thanks{First author partially supported by funds of the University of Salento, GNSAGA
and MIUR (PRIN)}
\date{}

\begin{abstract}
We consider three- and four-dimensional pseudo-Rieman\-nian generalized
symmetric spaces, whose invariant metrics were explicitly described in
\cite{CK}. While four-dimensional pseudo-Riemannian generalized symmetric
spaces of types $A$, $C$ and $D$ are algebraic Ricci solitons, the ones of
type $B$ are not so. The Ricci soliton equation for their metrics yields a
system of partial differential equations. Solving such system, we prove that
almost all the four-dimensional pseudo-Riemannian generalized symmetric spaces
of type $B$ are Ricci solitons. These examples show some deep differences
arising for the Ricci soliton equation between the Riemannian and the
pseudo-Riemannian cases, as any homogeneous Riemannian Ricci soliton is
algebraic \cite{Jab}. We also investigate three-dimensional generalized
symmetric spaces of any signature and prove that they are Ricci solitons.

\end{abstract}
\maketitle

\section{Introduction}

\setcounter{equation}{0}

Generalized symmetric spaces are a natural generalization of symmetric
spaces.
%We may refer to the monograph \cite{K} for an overview of the topic.
Since their introduction, the geometry of generalized symmetric spaces has
been intensively studied by several authors. Finite order automorphisms of
semisimple Lie algebras and Riemannian manifolds with geodesic symmetries of
order $3$ were studied respectively in \cite{Ka} and \cite{Gr}. In \cite{K},
O. Kowalski undertook a study of generalized symmetric spaces without using
neither topological invariants nor advanced algebra. Homogeneous structures of
generalized symmetric Riemannian spaces were studied in \cite{GC}. S.~Terzi\'c
classified generalized symmetric spaces defined as quotients of compact simple
Lie groups, describing explicitly their real cohomology algebras \cite{T1} and
calculating their real Pontryagin characteristic classes \cite{T2}. Formality
of all generalized symmetric spaces was proved by D.~Kotschick and S.~Terzi\'c
in \cite{KT}.

\v{C}ern\'{y} and Kowalski \cite{CK} completely classified pseudo-Riemannian
generalized symmetric spaces of dimension $\leq4$. In dimension $n=2$ they are
necessarily symmetric. In dimension $n=3$ the proper examples may be described
as $\mathbb{R}^{3}$ endowed with a special metric, with all possible
signatures. In dimension $n=4$, a proper generalized symmetric space may be
identified with $\mathbb{R}^{4}$ endowed with a special metric of four
different types, called in \cite{CK} types $A$, $B$, $C$ and $D$. The metrics
of type $A$ are either Riemannian or of neutral signature $(2,2)$; metrics of
type $B$ and $D$ always have signature $(2,2)$; for type $C$ (which was proved
in \cite{DV} to be indeed symmetric), the metric is Lorentzian. All these
spaces are reductive homogeneous.

Many aspects of the geometry of four-dimensional pseudo-Riemannian generalized
symmetric spaces have been investigated: K\"ahler and para-K\"ahler structures
\cite{Calvaruso}, harmonicity properties of vector fields \cite{Cejm},
curvature properties \cite{CD}, parallel degenerate distributions \cite{CZ0},
homogeneous geodesics \cite{DM}, parallel hypersurfaces \cite{DV}.
Four-dimensional pseudo-Riemannian generalized symmetric spaces of type $A$,
$C$ and $D$ are algebraic Ricci solitons, whereas those of type $B$ are never
algebraic Ricci solitons \cite{BO}.

A \emph{Ricci soliton} is a pseudo-Riemannian manifold $(M,g)$, together with
a smooth vector field $X$, such that
\begin{equation}
\mathcal{L}_{X}g+\varrho=\lambda g, \label{RicciSoliton}%
\end{equation}
where $\mathcal{L}_{X}$ and $\varrho$ respectively denote the Lie derivative
in the direction of $X$ and the Ricci tensor and $\lambda$ is a real number. A
Ricci soliton is said to be either \emph{shrinking}, \emph{steady} or
\emph{expanding}, according to whether $\lambda>0$, $\lambda=0$ or $\lambda
<0$, respectively. When $X$ is the gradient of some smooth function $f\colon
M\rightarrow\mathbb{R}$, the metric $g$ is said to be a \emph{gradient Ricci
soliton}.

Complete Ricci solitons with respect to some complete vector field correspond
to the self-similar solutions of the Ricci flow they generate. As such, they
play an essential role in understanding the singularities of the Ricci flow.
We may refer to the recent survey \cite{Cao} for more information and further
references on Ricci solitons. Pseudo-Riemannian Ricci solitons have been
recently studied by several authors, some examples may be found in
\cite{BO}-\cite{BGGG},\cite{CalNA},\cite{CD2}-\cite{CF1},\cite{CZ},\cite{On}
and references therein.

If we start from the explicit description of a pseudo-Riemannian metric $g$
with respect to some (local) coordinates, the Ricci soliton equation
(\ref{RicciSoliton}) leads to a system of partial differential equations.
However, in the case of homogeneous pseudo-Riemannian metrics, the first
approach in the study of the Ricci soliton equation (\ref{RicciSoliton}) is
usually algebraic. A \emph{homogeneous Ricci soliton} is a homogeneous space
$M=G/H$, together with a $G$-invariant metric $g$, for which equation
(\ref{RicciSoliton}) holds. An \emph{invariant Ricci soliton} is a homogeneous
one, such that equation (\ref{RicciSoliton}) holds for a $G$-invariant vector field.

Algebraic Ricci solitons were introduced by Lauret \cite{Lau} for Riemannian
manifolds and successively extended to pseudo-Riemannian settings \cite{On}.
Consider a homogeneous (reductive) pseudo-Riemannian manifold $(M=G/H,g)$ and
the corresponding reductive decomposition $\mathfrak{g}=\mathfrak{m}%
\oplus\mathfrak{h}$ of the Lie algebra $\mathfrak{g}$ of $G$. The metric $g$
is said to be an \emph{algebraic Ricci soliton} if there exists some
derivation $D\in\mathrm{Der}(\mathfrak{g})$, such that
\begin{equation}
Ric=c\,\mathrm{Id}+pr\circ D, \label{ARS}%
\end{equation}
where $Ric$ denotes the Ricci operator of $\mathfrak{m}$, $pr\colon
\mathfrak{g}\rightarrow\mathfrak{m}$ the projection and $c$ is a real number.
An algebraic Ricci soliton on a solvable Lie group is called a
\emph{solvsoliton}.

Any algebraic Ricci soliton metric $g$ is also a Ricci soliton (\cite{Lau}%
,\cite{On}), satisfying \eqref{RicciSoliton} with $\lambda=c$. We emphasize
the fact that in Riemannian settings this algebraic approach is essentially
all that is needed. In fact, homogeneous \emph{Riemannian} Ricci solitons are
indeed algebraic, with respect to some suitable transitive group $G$ of
isometries \cite{Jab}. Therefore, the algebraic approach to the study of Ricci
solitons exhausts all possibilities in the Riemannian case. On the other hand,
in pseudo-Riemannian settings, neither algebraic nor invariant Ricci solitons
are generally enough to determine all homogeneous Ricci solitons.
%Explicit examples of non-algebraic, non-invariant Lorentzian Ricci solitons may be found in \cite{CalNA}.
Finally, we observe that also in the case of an algebraic Ricci soliton
metric, this information does not codify all possible remarkable properties of
this solution. In fact, given an algebraic Ricci soliton, we do not know
whether it is also invariant or a gradient one.

In this paper we shall complete the study of the Ricci soliton equation
\eqref{RicciSoliton} on low-dimensional pseudo-Riemannian generalized
symmetric spaces, considering the three-dimensional examples and the
four-dimensional pseudo-Riemannian generalized symmetric spaces of type $B$.
Three-dimensional examples, of any signature, turn out to be algebraic Ricci
solitons. Retrieving them as explicit solutions of equation
\eqref{RicciSoliton}, we shall prove that they are
%neither invariant nor
not gradient. Four-dimensional generalized symmetric spaces of type $B$ will
provide a new family of examples of non-algebraic homogeneous
pseudo-Riemannian Ricci solitons.

{The paper is organized in the following way. In Section 2 we shall report the
description in a set of (global) coordinates of three-dimensional Lorentzian
generalized symmetric spaces and four-dimensional pseudo-Riem\-an\-nian
generalized symmetric spaces of type $B$, and we shall compute all the
curvature information with respect to the corresponding basis of coordinate
vector fields}.
%that we denote by $\partial_{1}=\partial/\partial x$, $\partial_{2}=\partial/\partial y$, $\partial_{3}=\partial/\partial u$, $\partial_{4}=\partial/\partial v$
In Sections 3 and 4 we shall introduce and solve the systems of PDEs that
express the Ricci soliton equation (\ref{RicciSoliton}) in these global coordinates.

\section{Generalized Symmetric spaces}

\setcounter{equation}{0}

Let $(M,g)$ be a pseudo-Riemannian manifold. An \emph{$s$-structure} on $M$ is
a family of isometries $\{s_{p}~|~p\in M\}$ (called \emph{symmetries}) of
$(M,g)$, such that each $s_{p}$ has $p$ as an isolated fixed point. The
\emph{order} of an $s$-structure is the least integer $k$ such that
$(s_{p})^{k}=\mathrm{Id}$ for all $p \in M$ (if such an integer does not
exist, the $s$-structure is said to be \emph{of infinite order}).

If, for any pair of points $p,q \in M$,

\begin{itemize}
\item[(i)] the mapping $(p,q)\to s_{p} (q)$ is smooth, and \vspace{4pt}

\item[(ii)] $s_{p} \circ s_{q}=s_{\tilde{q}} \circ s_{p}$, where $\tilde
{q}=s_{p} (q)$,
\end{itemize}

then the $s$-structure is said to be \emph{regular}. Each symmetric space
admits a regular $s$-structure, given by the family of its involutive geodesic
symmetries. More in general, a pseudo-Riemannian manifold admitting at least
one regular $s$-structure is called a \emph{generalized symmetric space}. The
infimum of all integers $k(\geq2)$ such that $(M,g)$ admits an $s$-structure
of order $k$, is called the \emph{order} of the generalized symmetric space.

\medskip Following \cite{CK}, any proper (that is, non-symmetric)
three-dimensional generalized symmetric space $(M,g)$ is of order $4$.
Moreover, it is given by the space $\mathbb{R}^{3}(x,y,t)$ with the
pseudo-Riemannian metric
\begin{equation}
\label{3D}g=\varepsilon\left(  e^{2t} dx^{2}+e^{-2t} dy^{2} \right)  + \mu
dt^{2} ,
\end{equation}
where $\varepsilon= \pm1$ and $\mu\neq0$ is a real constant. With respect to
the notations introduced in \cite{CK}, in equation \eqref{3D} we used $\mu$
instead of $\lambda$, which we reserved for the Ricci soliton equation
\eqref{RicciSoliton}. Depending on the values of $\varepsilon$ and $\mu$,
these metrics can attain any possible signature: $(3,0)$, $(0,3)$, $(2,1)$,
$(1,2)$.

\medskip As proved in \cite{CK}, a generalized symmetric space $(M,g)$ of type
$B$ is of order $3$ and is given by the space $\mathbb{R}^{4}(x,y,u,v)$ with
the pseudo-Riemannian metric
\begin{equation}
\label{TypeB}g=\mu\left(  dx^{2}+dy^{2}+dxdy\right)  +e^{-y}\left(
2dx+dy\right)  dv+e^{-x}\left(  dx+2dy\right)  du,
\end{equation}
where $\mu$ is a real constant. All these metrics have neutral signature
$(2,2)$.

\section{Ricci solitons on $3D$ generalized symmetric spaces}

\setcounter{equation}{0}

Starting from the description of the invariant metrics of three-dimensional
generalized symmetric spaces given in \eqref{3D}, we can explicitly calculate
their Levi-Civita connection and curvature. We use the coordinates
$(x_{1},x_{2},x_{3})=(x,y,t)$ and the corresponding basis of coordinate vector
fields $\left\{  \partial_{1},\partial_{2},\partial_{3}\right\}  =\left\{
\partial/\partial x, \partial/\partial y, \partial/\partial t \right\}  $.
Then, the Levi-Civita connection $\nabla$ is completely determined by the
following non-vanishing components:
\begin{equation}%
\begin{array}
[c]{ll}%
\nabla_{\partial_{1}}\partial_{1}=-\frac{\varepsilon}{\mu} e^{2x_{3}}
\partial_{3}, \smallskip & \nabla_{\partial_{1}}\partial_{3}= \partial_{1},
\smallskip\\[4pt]%
\nabla_{\partial_{2}}\partial_{2}=\frac{\varepsilon}{\mu} e^{-2x_{3}}
\partial_{3}, \smallskip & \nabla_{\partial_{2}}\partial_{3}=- \partial_{2}.
\end{array}
\label{Nabla3D}%
\end{equation}
We can now describe the Riemann-Christoffel curvature tensor $R$ of $(M,g)$
with respect to $\{\partial_{_{i}}\}$, computing $R(\partial_{i},\partial
_{j})=\nabla_{\partial_{i}}\nabla_{\partial_{j}}-\nabla_{\partial_{j}}%
\nabla_{\partial_{i}}$ for all indices $i,j$. Denoting by $R_{ij}$ the matrix
describing $R(\partial_{i},\partial_{j})$ with respect to the basis $\left\{
\partial_{1},\partial_{2},\partial_{3}\right\}  $ of coordinate vector fields,
we have
\begin{align*}
R_{12}  &  =\left(
\begin{array}
[c]{ccc}%
0 & \frac{\varepsilon}{\mu} e^{-2x_{3}} & 0\\
-\frac{\varepsilon}{\mu} e^{-2x_{3}} & 0 & 0\\
0 & 0 & 0\\
&  &
\end{array}
\right)  ,\\[6pt]
R_{13}  &  =\left(
\begin{array}
[c]{ccc}%
0 & 0 & -1\\
0 & 0 & 0\\
\frac{\varepsilon}{\mu} e^{2x_{3}} & 0 & 0
\end{array}
\right)  ,\qquad R_{23} =\left(
\begin{array}
[c]{ccc}%
0 & 0 & 0\\
0 & 0 & -1\\
0 & \frac{\varepsilon}{\mu} e^{2x_{3}} & 0
\end{array}
\right)  .
\end{align*}
Consequently, the Ricci tensor $\varrho(X,Y)=\mathrm{tr}(Z\mapsto R(Z,X)Y)$ is
then described with respect to the basis of coordinate vector fields by the
matrix
\begin{equation}
\varrho=\left(
\begin{array}
[c]{ccc}%
0 & 0 & 0\\[4pt]%
0 & 0 & 0\\[4pt]%
0 & 0 & -2
\end{array}
\right)  . \label{Ricci3D}%
\end{equation}
In particular, we observe that Ricci tensor is always degenerate.

With respect to the coordinate basis $\{\partial_{1},\partial_{2},\partial
_{3}\}$, let now $X=X^{i}\partial_{i}$ denote an arbitrary vector field on
$(M,g)$, where $X^{i}=X^{i}(x_{1},x_{2},x_{3})$, $i=1,2,3$, are arbitrary
smooth functions. We calculate $L_{X}g$ and obtain that the metric $g$,
together with the smooth vector field $X$, is a solution of the Ricci soliton
equation (\ref{RicciSoliton}) if and only if the following system of $6$ PDEs
is satisfied:
\begin{equation}
\left\{
\begin{array}
[c]{l}%
\varepsilon\,{\mathrm{e}^{2x_{{3}}}}\left(  2{\partial_{1}}X_{{1}}+2X_{{3}%
}-\lambda\right)  =0,\\[4pt]%
\varepsilon\,{\mathrm{e}^{2x_{{3}}}}{\partial_{{2}}}X_{{1}}+\varepsilon
\,\mathrm{e}^{-2x_{{3}}}\,{\partial_{{1}}}X_{{2}}=0,\\[4pt]%
\varepsilon\,{\mathrm{e}^{2x_{{3}}}}{\partial_{{3}}}X_{{1}}+\mu\,{\partial
_{{1}}}X_{{3}}=0,\\[4pt]%
\varepsilon\,{\mathrm{e}^{-2x_{{3}}}}\left(  2{\,{\partial_{{2}}}}X_{{2}%
}-2\,\,X_{{3}}-{\,\lambda}\right)  =0,\\[4pt]%
{\varepsilon{\mathrm{e}^{-2x_{{3}}}}\,{{\partial_{{3}}}}X_{{2}}}%
+\mu\,{{\partial_{{2}}}}X_{{3}}=0,\\[4pt]%
2\,\mu\,{{\partial_{{3}}}}X_{{3}}-2-\lambda\,\mu=0.
\end{array}
\right.  \label{3DSystem0}%
\end{equation}
We first integrate the last equation in \eqref{3DSystem0} and we get
\[
X_{{3}}\left(  x_{{1}},x_{{2}},x_{{3}}\right)  \,=\,{\frac{x_{{3}}}{\mu}%
}+\frac{1}{2}\,x_{{3}}\lambda+{F_{3}}\left(  x_{{1}},x_{{2}}\right)  ,
\]
where $F_{3}$ is an arbitrary smooth function. Substituting the above
expression of $X_{3}$ into the third and fifth equations of \eqref{3DSystem0},
they respectively become
\[%
\begin{array}
[c]{l}%
\varepsilon\,{\mathrm{e}^{2\,x_{{3}}}}{{\partial_{{3}}}}X_{{1}}+\mu
\,{{\partial_{{1}}}}{F_{3}}\left(  x_{{1}},x_{{2}}\right)  =0,\\[4pt]%
{\mathrm{e}^{-2\,x_{{3}}}}{{\partial_{{3}}}}X_{{2}}+\mu\,{{\partial_{{2}}}%
}{F_{3}}\left(  x_{{1}},x_{{2}}\right)  =0,
\end{array}
\]
which, integrated, yield
\[%
\begin{array}
[c]{l}%
X_{{1}}\left(  x_{{1}},x_{{2}},x_{{3}}\right)  \,=\frac{\mu}{2\varepsilon
}{\mathrm{e}^{-2\,x_{{3}}}}{{{{\partial_{{1}}}}{F_{3}}\left(  x_{{1}},x_{{2}%
}\right)  \mbox{}}}+{F_{1}}\left(  x_{{1}},x_{{2}}\right)  ,\\[8pt]%
X_{{2}}\left(  x_{{1}},x_{{2}},x_{{3}}\right)  =-\frac{\mu}{2\varepsilon
}\,{\mathrm{e}^{2\,x_{{3}}}}\,{{{{\partial_{{2}}}}{F_{3}}\left(  x_{{1}%
},x_{{2}}\right)  }}+{F_{2}}\left(  x_{{1}},x_{{2}}\right)  ,
\end{array}
\]
for some smooth functions $F_{1},F_{2}$. Substituting the above expressions of
$X_{1}$ and $X_{2}$ into the second equation of \eqref{3DSystem0}, it now
gives
\begin{equation}%
\begin{array}
[c]{l}%
\varepsilon\,{\mathrm{e}^{2\,x_{{3}}}}{{\partial_{{2}}}}{F_{1}}\left(  x_{{1}%
},x_{{2}}\right)  +\varepsilon\,{\mathrm{e}^{-2\,x_{{3}}}}{\partial_{{1}}%
}{F_{2}}\left(  x_{{1}},x_{{2}}\right)  =0.
\end{array}
\label{eq12}%
\end{equation}
Since the above equation \eqref{eq12} must hold for all values of $x_{3}$, it
yields ${{\partial_{{2}}}}{F_{1}}={\partial_{{1}}}{F_{2}}=0$, that is,
\[
{F_{1}}\left(  x_{{1}},x_{{2}}\right)  =G_{1}(x_{1}),\qquad{F_{2}}\left(
x_{{1}},x_{{2}}\right)  =G_{2}(x_{2}),
\]
for some smooth functions $G_{1},G_{2}$. System \eqref{3DSystem0} then reduces
to its first and fourth equations. In particular, the first equation in
\eqref{3DSystem0} now reads
\begin{equation}%
\begin{array}
[c]{l}%
\varepsilon{\mathrm{e}^{2\,x_{{3}}}}\left(  2\,\,{\partial_{1}G_{1}}(x_{{1}%
}){}+\left(  \frac{2}{\mu}+\lambda\right)  \,x_{3}\,+2\,{F_{3}}\left(  x_{{1}%
},x_{{2}}\right)  \mbox{}-\varepsilon\,\lambda\right) \\[4pt]%
+{{\mu}\,{{\partial_{11}^{2}}}}{F_{3}}\left(  x_{{1}},x_{{2}}\right)  =0.
\end{array}
\label{newsys}%
\end{equation}
Since the above equation \eqref{newsys} holds for all values of $x_{3}$, in
particular it implies that ${\partial_{11}^{2}}{F_{3}}=0$, from which, by
integration, we obtain
\[
F_{3}(x_{1},x_{2})=P_{3}(x_{2})x_{1}+Q_{3}(x_{2}),
\]
for some smooth functions $P_{3},Q_{3}$. We substitute into \eqref{newsys} and
it becomes
\[
{{\varepsilon\,{\mathrm{e}^{2\,x_{{3}}}}\left(  2\,{\partial_{{1}}}{G_{1}%
}\left(  x_{{1}}\right)  +\left(  \tfrac{2}{\mu}+\lambda\right)  x_{{3}%
}+2\,{P_{3}}\left(  x_{{2}}\right)  x_{{1}}+2\,{Q_{3}}\left(  x_{{2}}\right)
-\lambda\right)  =0}}.
\]
Again, this equation must hold for all values of $x_{3}$ and $x_{2}$.
Therefore, it implies at once that $\lambda=-\frac{2}{\mu}$ and that
$P_{3}(x_{2})=B_{3}$, $Q_{3}(x_{2})=A_{3}$, for some real constants
$A_{3},B_{3}$. Thus, the above equation now reduces to
\[
2\,{{\varepsilon\,{\mathrm{e}^{2\,x_{{3}}}}\left(  {\partial_{1}G_{1}}\left(
x_{{1}}\right)  +B_{{3}}x_{{1}}+A_{{3}}+\frac{1}{\mu}\right)  =0}}%
\]
and so, integrating we find
\[
{G_{1}}\left(  x_{{1}}\right)  \,=\,-\frac{1}{2}\,B_{{3}}x_{1}^{2}-\left(
A_{{3}}+{\frac{{{1}}}{\mu}}\right)  x_{1}+{A_{1}},
\]
where $A_{1}$ is a real constant. We are then left with the fourth equation in
\eqref{3DSystem0}, which now reads
\[
\,2\,{{\varepsilon\,{\mathrm{e}^{-2\,x_{{3}}}}\left(  \partial_{{2}}{G_{2}%
}(x_{{2}})-B_{{3}}x_{{1}}-A_{{3}}+\frac{1}{\mu}\right)  =0.}}%
\]
Since the above equation must hold for all values of $x_{1}$ and $x_{2}$, we
then necessarily have $B_{3}=0$ and, by integration,
\[
{G_{2}}(x_{2})\,=\left(  A_{3}-{\frac{1}{\mu}}\right)  x_{{2}}+{A_{2}},
\]
where $A_{2}$ is a real constant. All equations in system \eqref{3DSystem0}
are now satisfied. Therefore, the arbitrary invariant metric $g$ of a
three-dimensional generalized symmetric space is a Ricci soliton. More
precisely, replacing the explicit expressions of the functions we found above
into $X_{1},X_{2},X_{3}$, we conclude that, with respect to coordinates
$(x_{1},x_{2},x_{3})=(x,y,t)$, the Ricci soliton equation \eqref{RicciSoliton}
holds for the metric $g$ together with the vector field
\begin{equation}
X=\left(  A_{1}-\left(  A_{3}+\frac{1}{\mu}\right)  x_{1}\right)  \partial
_{1}+\left(  A_{2}+\left(  A_{3}-\frac{1}{\mu}\right)  x_{2}\right)
\partial_{2}+A_{3}\partial_{3}, \label{X3D}%
\end{equation}
where $A_{1},A_{2},A_{3}$ are real constants.

\smallskip In order to check the above result, we compute $\mathcal{L}_{X}g$
for the vector field $X$ described in \eqref{X3D} and compare it with
$-\varrho-\frac{2}{\mu}g$, which can be calculated at once from \eqref{3D} and
\eqref{Ricci3D}. We explicitly find, with respect to the basis $\left\{
\partial_{1}, \partial_{2}, \partial_{3} \right\}  $:%
\[
\mathcal{L}_{X}g=\left(
\begin{array}
[c]{ccc}%
-\frac{2\varepsilon}{\mu}e^{2x_{3}} & 0 & 0\\[4pt]%
0 & -\frac{2\varepsilon}{\mu}e^{-2x_{3}} & 0\\[4pt]%
0 & 0 & 0
\end{array}
\right)  =-\varrho-\frac{2}{\mu}g.
\]

\smallskip We now prove that the above Ricci soliton is not a gradient one,
that is, there no exists a smooth function $f(x_{1},x_{2},x_{3})$, such that
$X=\mathrm{grad}(f)=\sum_{i,j}g^{ij}\frac{\partial f}{\partial x_{i}}\partial
x_{j}$. In fact, suppose that such a function exists. Then, by \eqref{X3D} we
have
\begin{equation}%
\begin{array}
[c]{l}%
\varepsilon\,{\mathrm{e}^{-2\,x_{{3}}}}{{\partial_{{1}}}}f={A_{1}-{\left(
A_{{3}}+\frac{1}{\mu}\right)  \,x_{{1}}}},\\[8pt]%
\varepsilon\,{\mathrm{e}^{2\,x_{{3}}}}{{{\partial_{{2}}}}f}=A_{2}+\left(
A_{{3}}-\frac{1}{\mu}\right)  x_{{2}},\\[4pt]%
\frac{1}{\mu}{{\partial_{{3}}}}f=A_{{3}}.
\end{array}
\label{sysgrad3D}%
\end{equation}
We integrate the third equation of \eqref{sysgrad3D}, obtaining
\[
f(x_{1},x_{2},x_{3})=\mu A_{{3}}\,x_{{3}}+{W}(x_{{1}},x_{{2}}),
\]
where $W$ is a smooth function. We substitute into the second equation of
\eqref{sysgrad3D} and we get
\[
\varepsilon\,{\mathrm{e}^{2\,x_{{3}}}}{{\partial_{{2}}}}{W}\left(  x_{{1}%
},x_{{2}}\right)  =A_{2}+\left(  A_{{3}}-\frac{1}{\mu}\right)  x_{{2}}.
\]
Since the above equation must hold for all values of $x_{3}$, it implies that
${{\partial_{{2}}}}{W}\left(  x_{{1}},x_{{2}}\right)  =0$ and $A_{2}+\left(
A_{{3}}-\frac{1}{\mu}\right)  x_{{2}}=0$, for all $x_{2}$. Therefore,
$A_{2}=0$ and $A_{3}=\frac{1}{\mu}$) and integrating ${{\partial_{{2}}}}{W}=0$
we have
\[
W(x_{1},x_{2})=H(x_{1}),
\]
for a smooth function $H$. Finally, we substitute the above into the first
equation of \eqref{sysgrad3D}, which now reads
\[
\varepsilon{{\mathrm{e}^{-2\,x_{{3}}}}{\partial_{{1}}}{H}(x_{{1}})=A_{1}%
-\frac{2}{\mu}\,x_{{1}}}.
\]
Again, the above equation must be satisfied for all values of $x_{3}$. So,
${\partial_{{1}}}{H}(x_{{1}})=0$ (whence, $H(x_{1})=R_{1}$ is a real constant)
and we are left with
\[
A_{1}-\frac{2}{\mu}\,x_{{1}}=0,
\]
for all values of $x_{1}$, which cannot occur. Therefore, $g$ is never a
gradient Ricci soliton. This fact is not surprizing, since the existence of a
gradient non-steady Ricci soliton has some strong consequences on the geometry
of the manifold \cite{BGGG}.

The above results are summarized in the following.

\begin{theorem}
Every three-dimensional proper generalized symmetric space is a Ricci soliton.
More precisely, the arbitrary metric $g$ described by \eqref{3D} for
$\varepsilon=\pm1$ and $\mu\neq0$, satisfies the Ricci soliton equation
\eqref{RicciSoliton} taking
\[
%\begin{equation}\label{X3D}
X=\left(  A_{1}-\left(  A_{3}+\frac{1}{\mu}\right)  x\right)  \frac{\partial
}{\partial x}+\left(  A_{2}+\left(  A_{3}-\frac{1}{\mu}\right)  y\right)
\frac{\partial}{\partial y}+A_{3}\frac{\partial}{\partial t},
\]
%\end{equation}
where $A_{1},A_{2},A_{3}$ are real constants. This Ricci soliton is not a
gradient one.
\end{theorem}

\section{Ricci solitons on $4D$ generalized symmetric spaces of type $B$}

\setcounter{equation}{0}

We described the general pseudo-Riemannian metric of a four-dimensional
generalized symmetric space of type $B$ in \eqref{TypeB}. Proceeding similarly
to the three-dimensional case treated in the previous section, we first
compute its Levi-Civita connection and curvature. We consider global
coordinates $(x_{1},x_{2},x_{3},x_{4})=(x,y,u,v)$. With respect to the basis
of coordinate vector fields {\color{blue} }$\left\{  \partial_{1},\partial
_{2},\partial_{3},\partial_{4}\right\}  $, where $\partial_{i}=\partial
/\partial x_{i}${, }$1\leq i\leq4${,} the Levi-Civita connection $\nabla$ is
completely determined by the following non-vanishing components:
\begin{equation}
\left\{
\begin{array}
[c]{l}%
\nabla_{\partial_{1}}\partial_{1}=\frac{1}{3}\left(  \partial_{1}%
-2\partial_{2}\right)  +\frac{1}{3}\mu\left(  2{{\mathrm{e}}}^{x_{1}}%
\partial_{3}-{{\mathrm{e}}}^{x_{2}}\partial_{4}\right)  ,\smallskip\\
\nabla_{\partial_{1}}\partial_{2}=-\frac{1}{3}\left(  \partial_{1}%
+\partial_{2}\right)  +\frac{1}{3}\mu\left(  {{\mathrm{e}}}^{x_{1}}%
\partial_{3}+{{\mathrm{e}}}^{x_{2}}\partial_{4}\right)  ,\smallskip\\
\nabla_{\partial_{1}}\partial_{3}=-\frac{1}{3}{{\mathrm{e}}}^{-x_{1}}\left(
2{{\mathrm{e}}}^{x_{1}}\partial_{3}-{{\mathrm{e}}}^{x_{2}}\partial_{4}\right)
,\smallskip\\
\nabla_{\partial_{1}}\partial_{4}=\frac{1}{3}{{\mathrm{e}}}^{-x_{2}}\left(
2{{\mathrm{e}}}^{x_{1}}\partial_{3}-{{\mathrm{e}}}^{x_{2}}\partial_{4}\right)
,\smallskip\\
\nabla_{\partial_{2}}\partial_{2}=-\frac{1}{3}\left(  2\partial_{1}%
-\partial_{2}\right)  -\frac{1}{3}\mu\left(  {{\mathrm{e}}}^{x_{1}}%
\partial_{3}-2{{\mathrm{e}}}^{x_{2}}\partial_{4}\right)  ,\smallskip\\
\nabla_{\partial_{2}}\partial_{3}=-\frac{1}{3}{{\mathrm{e}}}^{-x_{1}}\left(
{{\mathrm{e}}}^{x_{1}}\partial_{3}-2{{\mathrm{e}}}^{x_{2}}\partial_{4}\right)
,\smallskip\\
\nabla_{\partial_{2}}\partial_{4}=\frac{1}{3}{{\mathrm{e}}}^{-x_{2}}\left(
{{\mathrm{e}}}^{x_{1}}\partial_{3}-2{{\mathrm{e}}}^{x_{2}}\partial_{4}\right)
.
\end{array}
\right.  \label{Nabla}%
\end{equation}

Consequently, we can then determine the Riemann-Christoffel curvature tensor
$R$ of $g$ with respect to $\{\partial_{_{i}}\}$. Denoting by $R_{ij}$ the
matrix describing $R(\partial_{i},\partial_{j})$ with respect to the basis
$\left\{  \partial_{1},\partial_{2},\partial_{3},\partial_{4}\right\}  $, of
coordinate vector fields, we have $R_{14}=R_{23}=R_{34}=0$ and
\begin{align*}
R_{12}  &  =\frac{1}{3}\left(
\begin{array}
[c]{cccc}%
-1 & -2 & 0 & 0\\
2 & 1 & 0 & 0\\
0 & 0 & 1 & 2e^{x_{1}-x_{2}}\\
0 & 0 & -2e^{-x_{1}+x_{2}} & -1
\end{array}
\right)  ,\\
R_{13}  &  =\frac{1}{3}\left(
\begin{array}
[c]{cccc}%
0 & 0 & 0 & 0\\
0 & 0 & 0 & 0\\
2 & 1 & 0 & 0\\
-e^{-x_{1}+x_{2}} & -2e^{-x_{1}+x_{2}} & 0 & 0
\end{array}
\right)  ,\\
R_{24}  &  =\frac{1}{3}\left(
\begin{array}
[c]{cccc}%
0 & 0 & 0 & 0\\
0 & 0 & 0 & 0\\
-2e^{x_{1}-x_{2}} & -e^{x_{1}-x_{2}} & 0 & 0\\
1 & 2 & 0 & 0
\end{array}
\right)  .
\end{align*}
{The Ricci tensor, }$\varrho(X,Y)=\mathrm{tr}(Z\mapsto R(Z,X)Y)$, is then
described with respect to the basis of coordinate vector fields by the matrix
\begin{equation}
\varrho=\left(
\begin{array}
[c]{cccc}%
-\frac{4}{3} & -\frac{2}{3} & 0 & 0\\[4pt]%
-\frac{2}{3} & -\frac{4}{3} & 0 & 0\\[4pt]%
0 & 0 & 0 & 0\\[4pt]%
0 & 0 & 0 & 0
\end{array}
\right)  \label{Ricci}%
\end{equation}
and so, is always degenerate.

Let now $X=X^{i}\partial_{i}$ denote an arbitrary vector field on $(M,g)$,
where $X^{i}=X^{i}(x_{1},x_{2},x_{3},x_{4})$, $1\leq i\leq4$, are arbitrary
smooth functions. We compute $L_{X}g$ and find that the metric $g$, together
with the smooth vector field $X$, is a solution of the Ricci soliton equation
(\ref{RicciSoliton}) if and only if the following system of $10$ PDEs is
satisfied:
\begin{equation}
\left\{
\begin{array}
[c]{crl}%
\!\!\!\text{(i)}\!\! & \mu\partial_{1}\left(  2X^{1}\!+\!X^{2}\right)
+e^{-x_{1}}\partial_{1}X^{3}+2e^{-x_{2}}\partial_{1}X^{4}\!-\!\tfrac{4}%
{3}\!-\!\lambda\mu & \!\!=0,\medskip\\
\!\!\!\text{(ii)}\!\! & \multicolumn{1}{l}{\mu\partial_{1}\left(
X^{1}\!+\!2X^{2}\right)  +\mu\partial_{2}\left(  2X^{1}\!+\!X^{2}\right) } &
\\
& +e^{-x_{2}}\left(  \partial_{1}+2\partial_{2}\right)  X^{4}+e^{-x_{1}%
}\left(  2\partial_{1}+\partial_{2}\right)  X^{3}\!-\!\tfrac{4}{3}%
\!-\!\lambda\mu & \!\!=0,\medskip\\
\!\!\!\text{(iii)}\!\! & \multicolumn{1}{l}{e^{-x_{1}}\partial_{1}\left(
X^{1}\!+\!2X^{2}\right)  +\mu\partial_{3}\left(  2X^{1}\!+\!X^{2}\right) } &
\\
& +\partial_{3}\left(  e^{-x_{1}}X^{3}+2e^{-x_{2}}X^{4}\right)  -e^{-x_{1}%
}X^{1}\!-\!\lambda e^{-x_{1}} & \!\!=0,\medskip\\
\!\!\!\text{(iv)}\!\! & \multicolumn{1}{l}{e^{-x_{2}}\partial_{1}\left(
2X^{1}\!+\!X^{2}\right)  +\mu\partial_{4}\left(  2X^{1}\!+\!X^{2}\right) } &
\\
& +\partial_{4}\left(  e^{-x_{1}}X^{3}+2e^{-x_{2}}X^{4}\right)
\!-\!2e^{-x_{2}}\left(  X^{2}\!+\!\lambda\right)  & \!\!=0,\medskip\\
\!\!\!\text{(v)}\!\! & \mu\partial_{2}\left(  X^{1}\!+\!2X^{2}\right)
+2e^{-x_{1}}\partial_{2}X^{3}+e^{-x_{2}}\partial_{2}X^{4}\!-\!\tfrac{4}%
{3}\!-\!\lambda\mu & \!\!=0,\medskip\\
\!\!\!\text{(vi)}\!\! & \multicolumn{1}{l}{e^{-x_{1}}\partial_{2}\left(
X^{1}\!+\!2X^{2}\right)  +\mu\partial_{3}\left(  X^{1}\!+\!2X^{2}\right) } &
\\
& +\partial_{3}\left(  2e^{-x_{1}}X^{3}+e^{-x_{2}}X^{4}\right)
\!-\!2e^{-x_{1}}\left(  X^{1}\!+\!\lambda\right)  & \!\!=0,\medskip\\
\!\!\!\text{(vii)}\!\! & \multicolumn{1}{l}{e^{-x_{2}}\partial_{2}\left(
2X^{1}\!+\!X^{2}\right)  +\mu\partial_{4}\left(  X^{1}\!+\!2X^{2}\right) } &
\\
& +\partial_{4}\left(  2e^{-x_{1}}X^{3}+e^{-x_{2}}X^{4}\right)  \!-\!e^{-x_{2}%
}\left(  X^{2}\!+\!\lambda\right)  & \!\!=0,\medskip\\
\!\!\!\text{(viii)}\!\! & e^{-x_{1}}\partial_{3}\left(  X^{1}\!+\!2X^{2}%
\right)  & \!\!=0,\medskip\\
\!\!\!\text{(ix)}\!\! & e^{-x_{2}}\partial_{3}\left(  2X^{1}\!+\!X^{2}\right)
+e^{-x_{1}}\partial_{4}\left(  X^{1}\!+\!2X^{2}\right)  & \!\!=0,\medskip\\
\!\!\!\text{(x)}\!\! & e^{-x_{2}}\partial_{4}\left(  2X^{1}\!+\!X^{2}\right)
& \!\!=0.
\end{array}
\right.  \label{System0}%
\end{equation}
Integrating equations (viii) and (x) in (\ref{System0}) we find that there
exist some smooth functions $Y^{1}\left(  x_{1},x_{2},x_{4}\right)  $ and
$=Y^{2}\left(  x_{1},x_{2},x_{3}\right)  $, such that
\begin{align}
X^{1}\left(  x_{1},x_{2},x_{3},x_{4}\right)  +2X^{2}\left(  x_{1},x_{2}%
,x_{3},x_{4}\right)   &  =Y^{1}\left(  x_{1},x_{2},x_{4}\right)  ,\label{1}\\
2X^{1}\left(  x_{1},x_{2},x_{3},x_{4}\right)  +X^{2}\left(  x_{1},x_{2}%
,x_{3},x_{4}\right)   &  =Y^{2}\left(  x_{1},x_{2},x_{3}\right)  . \label{2}%
\end{align}
Observe that $\partial_{3}Y^{2}$ does not depend on $x_{4}$, while
$\partial_{4}Y^{1}$ does not depend on $x_{3}$. Consequently, taking into
account (\ref{1}) and (\ref{2}), equation (ix) in (\ref{System0}) implies that
there exists some function $A(x_{1},x_{2})$, such that
\[
\tfrac{1}{2}e^{-x_{1}}\partial_{4}Y^{1}\left(  x_{1},x_{2},x_{4}\right)
=-\tfrac{1}{2}e^{-x_{2}}\partial_{3}Y^{2}\left(  x_{1},x_{2},x_{3}\right)
=A(x_{1},x_{2}).
\]
By integrating the above equations with respect to $Y^{1}(x_{1},x_{2},x_{3})$
and $Y^{2}(x_{1},x_{2},x_{4})$ respectively, we obtain%
\begin{align}
Y^{1}\left(  x_{1},x_{2},x_{4}\right)   &  =2e^{x_{1}}x_{4}A(x_{1}%
,x_{2})+B^{1}(x_{1},x_{2}),\label{3}\\
Y^{2}\left(  x_{1},x_{2},x_{3}\right)   &  =-2e^{x_{2}}x_{3}A(x_{1}%
,x_{2})+B^{2}(x_{1},x_{2}), \label{4}%
\end{align}
for some smooth functions $B^{i}(x_{1},x_{2})$, $i=1,2$. By (\ref{1}%
)-(\ref{4}), we now have
\begin{align}
X^{1}  &  =-\tfrac{2}{3}\left(  e^{x_{1}}x_{4}+2e^{x_{2}}x_{3}\right)
A(x_{1},x_{2})-\tfrac{1}{3}B^{1}(x_{1},x_{2})+\tfrac{2}{3}B^{2}(x_{1}%
,x_{2}),\label{X1bis}\\
X^{2}  &  =\tfrac{2}{3}\left(  2e^{x_{1}}x_{4}+e^{x_{2}}x_{3}\right)
A(x_{1},x_{2})+\tfrac{2}{3}B^{1}(x_{1},x_{2})-\tfrac{1}{3}B^{2}(x_{1},x_{2}).
\label{X2bis}%
\end{align}
Using (\ref{X1bis}) and (\ref{X2bis}), equation (iii) in (\ref{System0})
becomes
\begin{align*}
0  &  =\tfrac{1}{2}e^{-x_{1}}\partial_{3}X_{3}+e^{-x_{2}}\partial_{3}%
X_{4}+\partial_{1}A(x_{1},x_{2})x_{4}+\tfrac{1}{2}e^{-x_{1}}\partial_{1}%
B^{1}(x_{1},x_{2})\\
&  +\left(  \tfrac{4}{3}x_{4}+\tfrac{2}{3}x_{3}e^{-x_{1}+x_{2}}-\mu e^{x_{2}%
}\right)  A(x_{1},x_{2})\\
&  +\tfrac{1}{6}e^{-x_{1}}\left(  B^{1}(x_{1},x_{2})-2B^{2}(x_{1}%
,x_{2})\right)  -\tfrac{1}{2}\lambda e^{-x_{1}}.
\end{align*}
We integrate the above equation with respect to $X^{4}(x_{1},x_{2},x_{3}%
,x_{4})$ and we obtain
\begin{align}
X^{4}  &  =-\left(  \tfrac{1}{3}x_{3}e^{-x_{1}}+\tfrac{4}{3}x_{4}e^{-x_{2}%
}-\mu\right)  e^{2x_{2}}x_{3}A(x_{1},x_{2})\label{X4bis}\\
&  -e^{x_{2}}\partial_{1}A(x_{1},x_{2})x_{3}x_{4}-\tfrac{1}{2}\partial
_{1}B^{1}(x_{1},x_{2})e^{-x_{1}+x_{2}}x_{3}\nonumber\\
&  +\tfrac{1}{6}\left(  -B^{1}(x_{1},x_{2})+2B^{2}(x_{1},x_{2})\right)
e^{-x_{1}+x_{2}}x_{3}+\tfrac{1}{2}\lambda e^{-x_{1}+x_{2}}x_{3}\nonumber\\
&  -\tfrac{1}{2}e^{-x_{1}+x_{2}}X^{3}+Y^{4}(x_{1},x_{2},x_{4}),\nonumber
\end{align}
for some smooth function $Y^{4}(x_{1},x_{2},x_{4})$. Substituting from
(\ref{X1bis}), (\ref{X2bis}) and (\ref{X4bis}) into equation (vi) of
(\ref{System0}) it gives
\begin{align*}
0  &  =\tfrac{3}{4}e^{-x_{1}}\partial_{3}X^{3}+\tfrac{1}{4}e^{-x_{1}}\left(
2\partial_{2}B^{1}(x_{1},x_{2})-\partial_{1}B^{1}(x_{1},x_{2})\right) \\
&  -\tfrac{1}{2}\left(  \partial_{1}A(x_{1},x_{2})-2\partial_{2}A(x_{1}%
,x_{2})\right)  x_{4}\\
&  +\tfrac{1}{4}e^{-x_{1}}\left(  B^{1}(x_{1},x_{2})-2B^{2}(x_{1}%
,x_{2})\right) \\
&  +A(x_{1},x_{2})x_{3}e^{-x_{1}+x_{2}}+\tfrac{1}{2}\mu A(x_{1},x_{2}%
)e^{x_{2}}-\tfrac{3}{4}\lambda e^{-x_{1}}.
\end{align*}
We then integrate the above equation with respect to $X^{3}(x_{1},x_{2}%
,x_{3},x_{4})$ and we get
\begin{align}
X^{3}  &  =\tfrac{1}{3}\left(  \partial_{1}B^{1}(x_{1},x_{2})-2\partial
_{2}B^{1}(x_{1},x_{2})\right)  x_{3}\label{X3bis}\\
&  +\tfrac{2}{3}e^{x_{1}}\left(  \partial_{1}A(x_{1},x_{2})-2\partial
_{2}A(x_{1},x_{2})\right)  x_{3}x_{4}\nonumber\\
&  +\tfrac{1}{3}\left(  2B^{2}(x_{1},x_{2})-B^{1}(x_{1},x_{2})\right)
x_{3}-\tfrac{2}{3}A(x_{1},x_{2})e^{x_{2}}x_{3}^{2}\nonumber\\
&  -\tfrac{2}{3}A(x_{1},x_{2})e^{x_{2}}\mu e^{x_{1}}x_{3}+\lambda x_{3}%
+Y^{3}(x_{1},x_{2},x_{4}).\nonumber
\end{align}
Next, we substitute from (\ref{X1bis}), (\ref{X2bis}), (\ref{X4bis}) and
(\ref{X3bis}) into equation (iv) of (\ref{System0}) and we find
\begin{align*}
&  -2\left(  \partial_{1}A(x_{1},x_{2})+A(x_{1},x_{2})\right)  x_{3}\\
&  +e^{-x_{2}}\partial_{4}Y^{4}(x_{1},x_{2},x_{4})+\tfrac{1}{2}e^{-x_{2}%
}\partial_{1}B_{2}(x_{1},x_{2})\\
&  -\tfrac{4}{3}A(x_{1},x_{2})x_{4}e^{x_{1}-x_{2}}+\tfrac{1}{3}e^{-x_{2}%
}\left(  B_{2}(x_{1},x_{2})-2B^{1}(x_{1},x_{2})\right)  -\lambda e^{-x_{2}}=0.
\end{align*}
The above equation must be satisfied for the values of $x_{3}$. Therefore, it
yields
\begin{align}
&  \partial_{1}A(x_{1},x_{2})+A(x_{1},x_{2})=0,\label{5}\\[4pt]
&  e^{-x_{2}}\partial_{4}Y^{4}(x_{1},x_{2},x_{4})+\tfrac{1}{2}e^{-x_{2}%
}\partial_{1}B_{2}(x_{1},x_{2})-\tfrac{4}{3}A(x_{1},x_{2})x_{4}e^{x_{1}-x_{2}%
}\label{6}\\
&  +\tfrac{1}{3}e^{-x_{2}}\left(  B_{2}(x_{1},x_{2})-2B^{1}(x_{1}%
,x_{2})\right)  -\lambda e^{-x_{2}}=0.\nonumber
\end{align}
By integrating (\ref{5}) with respect to $A(x_{1},x_{2})$, we find
\begin{equation}
A(x_{1},x_{2})=C(x_{2})e^{-x_{1}}, \label{A}%
\end{equation}
for some smooth function $C(x_{2})$. We then substitute (\ref{A}) into
(\ref{6}) and integrate it with respect to $Y^{4}(x_{1},x_{2},x_{4})$,
obtaining
\begin{align}
Y^{4}(x_{1},x_{2},x_{4})  &  =\left(  -\tfrac{1}{2}\partial_{1}B^{2}%
(x_{1},x_{2})+\tfrac{1}{3}\left(  2B^{1}(x_{1},x_{2})-B^{2}(x_{1}%
,x_{2})\right)  \right. \label{Y4}\\
&  \left.  +\lambda\right)  x_{4}+\tfrac{2}{3}C(x_{2})x_{4}^{2}+B^{4}%
(x_{1},x_{2}),\nonumber
\end{align}
for some smooth function $B^{4}(x_{1},x_{2})$.

We now use (\ref{X1bis})-(\ref{X3bis}), (\ref{A}) and (\ref{Y4}) into equation
(vii) of (\ref{System0}), so that it yields
\begin{align*}
&  \tfrac{3}{4}e^{-x_{1}}\partial_{4}Y^{3}(x_{1},x_{2},x_{4})+\tfrac{1}%
{4}e^{-x_{2}}\left(  2\partial_{2}B_{2}(x_{1},x_{2})-\partial_{1}B^{2}%
(x_{1},x_{2})\right) \\
&  -2e^{-x_{1}}\left(  C^{\prime}(x_{2})+C(x_{2})\right)  x_{3}+\mu
C(x_{2})=0,
\end{align*}
or equivalently,
\begin{align*}
&  C^{\prime}(x_{2})+C(x_{2})=0,\\
&  3e^{-x_{1}}\partial_{4}Y^{3}(x_{1},x_{2},x_{4})+e^{-x_{2}}\left(
2\partial_{2}B_{2}(x_{1},x_{2})-\partial_{1}B^{2}(x_{1},x_{2})\right)
+\tfrac{4\mu}{3}C(x_{2})=0.
\end{align*}
By integrating the above equations with respect to $C(x_{2})$ and $Y^{3}%
(x_{1},x_{2},x_{4})$ respectively, we obtain $C(x_{2})=C_{2}e^{-x_{2}}$ and
\begin{align}
Y^{3}(x_{1},x_{2},x_{4})  &  =\tfrac{1}{3}e^{x_{1}-x_{2}}\left(  \partial
_{1}B^{2}(x_{1},x_{2})\!-\!2\partial_{2}B^{2}(x_{1},x_{2})\!-\!4\mu
C_{2}\right)  x_{4}\label{Y3}\\
&  \qquad\qquad\qquad\qquad\qquad\qquad\qquad+B^{3}(x_{1},x_{2}),\nonumber
\end{align}
for some smooth function $B^{3}(x_{1},x_{2})$.

Taking (\ref{X1bis})--(\ref{X3bis}) and (\ref{A})--(\ref{Y3}) into account,
equation (v) of (\ref{System0}) becomes
\begin{align}
&  -\tfrac{2}{3}C_{2}e^{-x_{2}}\left(  e^{-x_{2}}x_{4}+2e^{-x_{1}}%
x_{3}\right)  x_{4}-e^{-x_{1}}\left(  \partial_{2}^{2}B^{1}(x_{1}%
,x_{2})\right. \label{Eq22}\\
&  \left.  +\tfrac{2}{3}\left(  -2\partial_{2}B^{2}(x_{1},x_{2})+\partial
_{1}B^{1}(x_{1},x_{2})\right)  +\tfrac{1}{3}\partial_{2}B^{1}(x_{1}%
,x_{2})-\tfrac{4}{3}\mu C_{2}\right)  x_{3}\nonumber\\
&  -e^{-x_{2}}\left(  \partial_{2}^{2}B^{2}(x_{1},x_{2})-\tfrac{2}{3}%
\partial_{2}B^{1}(x_{1},x_{2})+\tfrac{2}{3}\partial_{1}B^{2}(x_{1}%
,x_{2})\right. \nonumber\\
&  \left.  -\partial_{2}B^{2}(x_{1},x_{2})-\tfrac{2}{3}\mu e^{x_{2}}%
C_{2}e^{-x_{2}}\right)  x_{4}+\tfrac{3}{2}e^{-x_{1}}\partial_{2}B^{3}%
(x_{1},x_{2})-\tfrac{4}{3}\nonumber\\
&  +e^{-x_{2}}\partial_{2}B^{4}(x_{1},x_{2})-\tfrac{1}{2}e^{-x_{1}}B^{3}%
(x_{1},x_{2})+\mu\left(  \partial_{2}B^{1}(x_{1},x_{2})-\lambda\right)  =0
.\nonumber
\end{align}
The above equation has to be satisfied for every value of $x_{3}$ and $x_{4}$.
Hence, we necessarily have $C_{2}=0$.

Thus, we are left with equations (i), (ii) and (v) of (\ref{System0}). Taking
into account (\ref{X1bis})--(\ref{X3bis}), (\ref{A})--(\ref{Y3}) and
$C(x_{2})=0$, these equations are rewritten as follows:
\begin{align}
&  e^{-x_{1}}\left(  \partial_{1}^{2}B^{1}(x_{1},x_{2})-\partial_{1}%
B^{1}(x_{1},x_{2})+\tfrac{2}{3}\partial_{2}B^{1}(x_{1},x_{2})\right.
\label{Ec11}\\
&  \left.  -\tfrac{2}{3}\partial_{1}B^{2}(x_{1},x_{2})\right)  x_{3}%
+e^{-x_{2}}\left(  \partial_{1}^{2}B^{2}(x_{1},x_{2})+\tfrac{2}{3}\partial
_{2}B^{2}(x_{1},x_{2})\right. \nonumber\\
&  \left.  +\tfrac{1}{3}\partial_{1}B^{2}(x_{1},x_{2})-\tfrac{4}{3}%
\partial_{1}B^{1}(x_{1},x_{2})\right)  x_{4}-2e^{-x_{2}}\partial_{1}%
B^{4}(x_{1},x_{2})\nonumber\\
&  -e^{-x_{1}}B^{3}(x_{1},x_{2})-\mu\partial_{1}B^{2}(x_{1},x_{2})+\tfrac
{4}{3}+\lambda\mu=0,\nonumber\\
&  e^{-x_{1}}\left(  \partial_{2}\partial_{1}B^{1}(x_{1},x_{2})+\tfrac{2}%
{3}\partial_{1}B^{1}(x_{1},x_{2})-\tfrac{2}{3}\partial_{1}B^{2}(x_{1}%
,x_{2})\right. \label{Ec12}\\
&  \left.  -\tfrac{1}{3}\partial_{2}B^{2}(x_{1},x_{2})\right)  x_{3}%
+e^{-x_{2}}\left(  \partial_{2}\partial_{1}B^{2}(x_{1},x_{2})+\tfrac{2}%
{3}\partial_{2}B^{2}(x_{1},x_{2})\right. \nonumber\\
&  \left.  -\tfrac{1}{3}\partial_{1}B^{1}(x_{1},x_{2})-\tfrac{2}{3}%
\partial_{2}B^{1}(x_{1},x_{2})\right)  x_{4}-\tfrac{1}{2}e^{-x_{2}}%
\partial_{1}B^{4}(x_{1},x_{2})\nonumber\\
&  -e^{-x_{2}}\partial_{2}B^{4}(x_{1},x_{2})-\tfrac{3}{4}e^{-x_{1}}%
\partial_{1}B^{3}(x_{1},x_{2})+\tfrac{1}{4}e^{-x_{1}}B^{3}(x_{1}%
,x_{2})\nonumber\\
&  +\tfrac{2}{3}-\tfrac{1}{2}\mu\left(  \partial_{1}B^{1}(x_{1},x_{2}%
)+\partial_{2}B^{2}(x_{1},x_{2})-\lambda\right)  =0 ,\nonumber\\
&  e^{-x_{1}}\left(  \partial_{2}^{2}B^{1}(x_{1},x_{2})+\tfrac{2}{3}%
\partial_{1}B^{1}(x_{1},x_{2})-\tfrac{4}{3}\partial_{2}B^{2}(x_{1}%
,x_{2})\right. \label{Ec22}\\
&  \left.  +\tfrac{1}{3}\partial_{2}B^{1}(x_{1},x_{2})\right)  x_{3}%
+e^{-x_{2}}\left(  \partial_{2}^{2}B^{2}(x_{1},x_{2})+\tfrac{2}{3}\partial
_{1}B^{2}(x_{1},x_{2})\right. \nonumber\\
&  \left.  -\tfrac{2}{3}\partial_{2}B^{1}(x_{1},x_{2})-\tfrac{2}{3}%
\partial_{2}B^{2}(x_{1},x_{2})\right)  x_{4}-\tfrac{3}{2}e^{-x_{1}}%
\partial_{2}B^{3}(x_{1},x_{2})\nonumber\\
&  -e^{-x_{2}}\partial_{2}B^{4}(x_{1},x_{2})+\tfrac{1}{2}e^{-x_{1}}B^{3}%
(x_{1},x_{2})+\tfrac{4}{3}\nonumber\\
&  -\mu\left(  \partial_{2}B^{1}(x_{1},x_{2})-\lambda\right)  =0 .\nonumber
\end{align}
Since equations (\ref{Ec11}), (\ref{Ec12}), (\ref{Ec22}) must hold for all
values of $x_{3}$ and $x_{4}$, they are equivalent to the following system of
equations:
\begin{equation}
\left\{
\begin{array}
[c]{cll}%
\text{(i)} & \partial_{1}^{2}B^{1}(x_{1},x_{2})-\partial_{1}B^{1}(x_{1}%
,x_{2})+\tfrac{2}{3}\partial_{2}B^{1}(x_{1},x_{2}) & \smallskip\\[3pt]%
\multicolumn{1}{r}{} & \multicolumn{1}{r}{-\tfrac{2}{3}\partial_{1}B^{2}%
(x_{1},x_{2})} & \!\!=0,\smallskip\\[6pt]%
\text{(ii)} & \partial_{2}\partial_{1}B^{1}(x_{1},x_{2})+\tfrac{2}{3}%
\partial_{1}B^{1}(x_{1},x_{2}) & \smallskip\\[3pt]%
\multicolumn{1}{r}{} & \multicolumn{1}{r}{-\tfrac{2}{3}\partial_{1}B^{2}%
(x_{1},x_{2})-\tfrac{1}{3}\partial_{2}B^{2}(x_{1},x_{2})} & \!\!=0,\smallskip
\\[6pt]%
\text{(iii)} & \partial_{2}^{2}B^{1}(x_{1},x_{2})+\tfrac{2}{3}\partial
_{1}B^{1}(x_{1},x_{2})+\tfrac{1}{3}\partial_{2}B^{1}(x_{1},x_{2}) &
\smallskip\\[3pt]%
\multicolumn{1}{r}{} & \multicolumn{1}{r}{-\tfrac{4}{3}\partial_{2}B^{2}%
(x_{1},x_{2})} & \!\!=0,\smallskip\\[6pt]%
\text{(iv)} & \partial_{1}^{2}B^{2}(x_{1},x_{2})+\tfrac{2}{3}\partial_{2}%
B^{2}(x_{1},x_{2})+\tfrac{1}{3}\partial_{1}B^{2}(x_{1},x_{2}) & \smallskip
\\[3pt]%
\multicolumn{1}{r}{} & \multicolumn{1}{r}{-\tfrac{4}{3}\partial_{1}B^{1}%
(x_{1},x_{2})} & \!\!=0,\smallskip\\[6pt]%
\text{(v)} & \partial_{2}\partial_{1}B^{2}(x_{1},x_{2})+\tfrac{2}{3}%
\partial_{2}B^{2}(x_{1},x_{2}) & \smallskip\\[3pt]%
\multicolumn{1}{r}{} & \multicolumn{1}{r}{-\tfrac{1}{3}\partial_{1}B^{1}%
(x_{1},x_{2})-\tfrac{2}{3}\partial_{2}B^{1}(x_{1},x_{2})} & \!\!=0,\smallskip
\\[6pt]%
\text{(vi)} & \partial_{2}^{2}B^{2}(x_{1},x_{2})+\tfrac{2}{3}\partial_{1}%
B^{2}(x_{1},x_{2})-\tfrac{2}{3}\partial_{2}B^{2}(x_{1},x_{2}) & \smallskip
\\[3pt]%
\multicolumn{1}{r}{} & \multicolumn{1}{r}{-\tfrac{2}{3}\partial_{2}B^{1}%
(x_{1},x_{2})} & \!\!=0,\smallskip\\[6pt]%
\text{(vii)} & -2e^{-x_{2}}\partial_{1}B^{4}(x_{1},x_{2})-e^{-x_{1}}%
B^{3}(x_{1},x_{2}) & \smallskip\\[3pt]%
\multicolumn{1}{r}{} & \multicolumn{1}{r}{+\tfrac{4}{3}-\mu\left(
\partial_{1}B^{2}(x_{1},x_{2})-\lambda\right) } & \!\!=0,\smallskip\\[6pt]%
\text{(viii)} & e^{-x_{2}}\left(  \tfrac{1}{2}\partial_{1}B^{4}(x_{1}%
,x_{2})+\partial_{2}B^{4}(x_{1},x_{2})\right)  & \smallskip\\[3pt]%
\multicolumn{1}{r}{} & \multicolumn{1}{r}{+e^{-x_{1}}\left(  \tfrac{3}%
{4}\partial_{1}B^{3}(x_{1},x_{2})-\tfrac{1}{4}B^{3}(x_{1},x_{2})\right)
-\tfrac{2}{3}} & \smallskip\\[3pt]%
\multicolumn{1}{r}{} & \multicolumn{1}{r}{+\tfrac{1}{2}\mu\left(  \partial
_{1}B^{1}(x_{1},x_{2})+\partial_{2}B^{2}(x_{1},x_{2})-\lambda\right) } &
\!\!=0,\smallskip\\[6pt]%
\text{(ix)} & e^{-x_{1}}\left(  \tfrac{3}{2}\partial_{2}B^{3}(x_{1}%
,x_{2})-\tfrac{1}{2}B^{3}(x_{1},x_{2})\right)  -\tfrac{4}{3} & \smallskip
\\[3pt]
& \multicolumn{1}{r}{+e^{-x_{2}}\partial_{2}B^{4}(x_{1},x_{2})+\mu\left(
\partial_{2}B^{1}(x_{1},x_{2})-\lambda\right) } & \!\!=0,
\end{array}
\right.  \label{System2}%
\end{equation}
{In discussing the solutions of system \eqref{System2}, we shall treat
separately the cases $\mu\neq0$ and $\mu=0$.}

\subsection{Case $\mu\neq0$}

{In this case we shall provide an explicit solution. In fact, suppose that
there exist some smooth functions }$P^{i}(x_{1})$, $Q^{i}(x_{2})$ {such that}
\begin{equation}
B^{i}(x_{1},x_{2})=P^{i}(x_{1})+Q^{i}(x_{2}),\quad i=1,2.
\label{SpecialSolution}%
\end{equation}
Taking into account (\ref{SpecialSolution}), from the equation (i) in
(\ref{System2}) we deduce that there exists some real constant $J^{1}$, such
that%
\[
\tfrac{2}{3}\frac{d}{dx_{2}}Q^{1}(x_{2})=\tfrac{2}{3}\frac{d}{dx_{1}}%
P^{2}(x_{1})-\frac{d^{2}}{dx^{2}}P^{1}(x_{1})+\frac{d}{dx_{1}}P^{1}%
(x_{1})=\tfrac{2}{3}J^{1}.
\]
Therefore,%
\[
Q^{1}(x_{2})=J^{1}x_{2}+H^{1},
\]
{where $H^{1}$ is a real constant. Replacing the above expressions of
$Q(x_{2})$ and of $B^{i}(x,x_{2})$, $i=1,2$,} equation (iii) of (\ref{System2}%
) becomes
\[
0=-\tfrac{2}{3}\frac{d}{dx_{1}}P^{1}(x_{1})+\tfrac{4}{3}\frac{d}{dx_{2}}%
Q^{2}(x_{2})-\tfrac{1}{3}J^{1},
\]
{which in particular yields} $\frac{d^{2}}{dx_{2}^{2}}Q^{2}(x_{2})=\frac
{d^{2}}{dx_{1}^{2}}P^{1}(x_{1})=0$. Hence, {integrating we get}
\[
Q^{2}(x_{2})=J^{2}x_{2}+H^{2}\text{ and }P^{1}(x_{1})=R^{1}x_{1}+S^{1},
\]
{for some real constants $J^{2}$, $H^{2}$, $R^{1}$ and }$S^{1}${. Then,
equations (i)-(vi) in (\ref{System2}) now give}
\[
\left\{
\begin{array}
[c]{l}%
3R^{1}-J^{1}+\frac{d}{dx_{1}}P^{2}(x_{1})=0,\\[6pt]%
2R^{1}-J^{2}-2\frac{d}{dx}P^{2}(x_{1})=0,\\[6pt]%
2R^{1}-4J^{2}+J^{1}=0,\\[6pt]%
3\frac{d^{2}}{dx_{1}^{2}}P^{2}(x_{1})+\frac{d}{dx_{1}}P^{2}(x_{1}%
)-4R^{1}+2J^{2}=0,\\[6pt]%
2J^{1}-2J^{2}+R^{1}=0,\\[6pt]%
J^{1}+3J^{2}-2\frac{d}{dx_{1}}P^{2}(x_{1})=0,
\end{array}
\right.
\]
{which easily yield {$R^{1}={{J^{2}=}}$ $J^{1}$}$=0$ and} $P^{2}(x_{1})=S^{2}$
for some real constant $S^{2}$. Therefore, {$B^{i}(x_{1},x_{2})=S^{i}+H^{i}$,
$i=1,2$, where $S^{i},H^{i}$ are real constants. Setting $S^{i}+H^{i}=W^{i}$,
we then have
\[
B^{i}(x_{1},x_{2})=W^{i},\qquad i=1,2,
\]
} and system (\ref{System2}) reduces to:%
\begin{equation}
\left\{
\begin{array}
[c]{l}%
2e^{-x_{2}}\partial_{1}B^{4}(x_{1},x_{2})+e^{-x_{1}}B^{3}(x_{1},x_{2}%
)-\lambda\mu-\tfrac{4}{3}=0,\smallskip\\[6pt]%
e^{-x_{1}}\left(  \tfrac{3}{2}\partial_{1}B^{3}(x_{1},x_{2})-\frac{1}{2}%
B^{3}(x_{1},x_{2})\right) \\[4pt]%
\multicolumn{1}{r}{\quad+e^{-x_{2}}\left(  \partial_{1}B^{4}(x_{1}%
,x_{2})+2\partial_{2}B^{4}(x_{1},x_{2})\right)  -\lambda\mu-\tfrac{4}%
{3}=0,\smallskip}\\[6pt]%
e^{-x_{1}}\left(  \frac{3}{2}\partial_{2}B^{3}(x_{1},x_{2})-\frac{1}{2}%
B^{3}(x_{1},x_{2})\right) \\[4pt]%
\multicolumn{1}{r}{+e^{-x_{2}}\partial_{2}B^{4}(x_{1},x_{2})-\lambda\mu
-\tfrac{4}{3}=0.}%
\end{array}
\right.  \label{System3}%
\end{equation}
From the first equation of \eqref{System3} we get
\[
B^{3}(x_{1},x_{2})=\left(  \lambda\mu+\tfrac{4}{3}\right)  e^{x_{1}}%
-2e^{x_{1}-x_{2}}\partial_{1}B^{4}(x_{1},x_{2}).
\]
Using the above expression for $B^{3}(x_{1},x_{2})$ and assuming $B^{4}%
(x_{1},x_{2})=P^{4}(x_{1})+Q^{4}(x_{2})$, from the second equation in
(\ref{System3}) we deduce
\[
Q^{4}(x_{2})=W^{4}\text{ and }P^{4}(x_{1})=R^{4}e^{-\frac{1}{3}x_{1}}+S^{4},
\]
where{ $W^{4}$, }$R^{4}$ and $S^{4}${ are a real constants. We then
substitute} the expressions for $B^{i}(x_{1},x_{2})$, $i=3,4$ in the last
equation of (\ref{System3}) and we obtain $R^{4}=0$ and $\lambda=-\tfrac
{4}{3\mu}$, where we clearly see that the condition $\mu\neq0$ is needed for
such a kind of solution.

{All equations in \eqref{System2} are then satisfied. We now replace the
functions we found by integration into $X^{i}(x_{1},x_{2},x_{3},x_{4})$. With
respect to global coordinates }$\left(  {x_{1},x_{2},x_{3},x_{4}}\right)  ${,
we then get a solution $X\left(  x_{1},x_{2},x_{3},x_{4}\right)
=X^{i}\partial_{i}$ of the Ricci soliton equation \eqref{RicciSoliton},
explicitly given by
\begin{equation}
\left\{
\begin{array}
[c]{l}%
X^{1}\left(  x_{1},x_{2},x_{3},x_{4}\right)  =-\tfrac{1}{3}W^{1}+\tfrac{2}%
{3}W^{2},\smallskip\\[4pt]%
X^{2}\left(  x_{1},x_{2},x_{3},x_{4}\right)  =\tfrac{2}{3}W^{1}-\tfrac{1}%
{3}W^{2},\smallskip\\[4pt]%
X^{3}\left(  x_{1},x_{2},x_{3},x_{4}\right)  =(-\tfrac{1}{3}W^{1}+\tfrac{2}%
{3}W^{2})x_{3}-\tfrac{4}{3\mu}x_{3},\smallskip\\[4pt]%
X^{4}\left(  x_{1},x_{2},x_{3},x_{4}\right)  =(\frac{2}{3}W^{1}-\frac{1}%
{3}W^{2})x_{4}-\tfrac{4}{3\mu}x_{4}+W^{3},\smallskip
\end{array}
\right.  \label{munot0}%
\end{equation}
for some real constants $W^{1},W^{2},W^{4}$ and }$W^{3}=S^{4}+W^{4}${. }

{As a check, by \eqref{munot0} and computing $\mathcal{L}_{X}g$, we explicitly
find, with respect to the basis $\left(  \partial_{i}\right)  _{1\leq i\leq4}%
$:}%
\[
\mathcal{L}_{X}g=\left(
\begin{array}
[c]{cccc}%
0 & 0 & -\frac{2}{3\mu}e^{-x_{1}} & -\frac{4}{3\mu}e^{-x_{2}}\\[4pt]%
0 & 0 & -\frac{4}{3\mu}e^{-x_{1}} & -\frac{2}{3\mu}\lambda e^{-x_{2}}\\[4pt]%
-\frac{2}{3\mu}e^{-x_{1}} & -\frac{4}{3\mu}e^{-x_{1}} & 0 & 0\\[4pt]%
-\frac{4}{3\mu}e^{-x_{2}} & -\frac{2}{3\mu}e^{-x_{2}} & 0 & 0
\end{array}
\right)  =-\varrho-\frac{4}{3\mu}g.
\]

Therefore, the arbitrary invariant metric $g$ for $\mu\neq0$ of a
four-dimensional generalized symmetric space of type $B$ is a Ricci soliton.
More precisely, replacing the explicit expressions of the functions we found
above into $X_{1}$, $X_{2}$, $X_{3}$, $X_{4}$, we conclude that, with respect
to coordinates $(x_{1},x_{2},x_{3},x_{4})=(x,y,u,v)$, the Ricci soliton
equation \eqref{RicciSoliton} holds for the metric $g$ with $\mu\neq0$,
together with the vector field
\begin{align}
X  &  =\tfrac{1}{3}\left(  -W^{1}+2W^{2}\right)  \partial_{1}+\tfrac{1}%
{3}\left(  2W^{1}-W^{2}\right)  \partial_{2}-\tfrac{1}{3}(W^{1}\!-\!2W^{2}%
)x_{3}\partial_{3}\label{X4D}\\
&  -\tfrac{4}{3\mu}x_{3}\partial_{3}+\tfrac{1}{3}\left(  \left(
2W^{1}\!-\!W^{2}\right)  x_{4}-\tfrac{4}{\mu}x_{4}+W^{3}\right)  \partial
_{4},\nonumber
\end{align}
where $W^{1}$ $W^{2}$ and $W^{3}$ are real constants.

\smallskip We now check that the above Ricci soliton is not a gradient one,
that is, there no exists a smooth function $f(x_{1},x_{2},x_{3},x_{4})$, such
that $X=\mathrm{grad}(f)=\sum_{i,j}g^{ij}\frac{\partial f}{\partial x_{i}%
}\partial_{j}$. In fact, suppose that such a function exists. Then, by
\eqref{X4D} we have%
\begin{equation}%
\begin{array}
[c]{l}%
-2{\mathrm{e}}^{x_{1}}\partial_{3}f+4{\mathrm{e}}^{x_{2}}\partial_{4}%
f+W^{1}-2W^{2}=0,\\[6pt]%
4{\mathrm{e}}^{x_{1}}\partial_{3}f-2{\mathrm{e}}^{x_{2}}\partial_{4}%
f-2W^{1}+W^{2}=0,\\[6pt]%
2{\mathrm{e}}^{x_{1}}\left(  -\partial_{1}f+2\partial_{2}f-\mu\left(
2{\mathrm{e}}^{x_{1}}\partial_{3}f-{\mathrm{e}}^{x_{2}}\partial_{4}f\right)
\right) \\[3pt]%
\multicolumn{1}{r}{+(W^{1}-2W^{2}+\tfrac{4}{\mu})x_{3}=0,}\\[6pt]%
2{\mathrm{e}}^{x_{2}}\left(  \partial_{2}f\!-\!2\partial_{1}f\!-\!\mu\left(
{\mathrm{e}}^{x_{1}}\partial_{3}f\!-\!2{\mathrm{e}}^{x_{2}}\partial
_{4}f\right)  \right) \\
\multicolumn{1}{r}{\!-\!\left(  2W^{1}\!-\!W^{2}\!-\!\tfrac{4}{\mu}\right)
x_{4}\!-\!3W^{3}=0.}%
\end{array}
\label{GradientRS4}%
\end{equation}
From the first equation of (\ref{GradientRS4}) we have%
\[
\partial_{3}f=2{\mathrm{e}}^{x_{2}-x_{1}}\partial_{4}f+{\mathrm{e}}^{-x_{1}%
}\left(  \tfrac{1}{2}W^{1}-W^{2}\right)  ,
\]
and substituting the above expression into the second equation of
(\ref{GradientRS4}), we get $\partial_{4}f=\tfrac{1}{2}{\mathrm{e}}^{-x_{2}%
}W^{2}$. Therefore, $\partial_{3}f=\tfrac{1}{2}{\mathrm{e}}^{-x_{1}}W^{1}$,
which by integration yields
\[
f(x_{1},x_{2},x_{3},x_{4})=\tfrac{1}{2}{\mathrm{e}}^{-x_{1}}W^{1}x_{{3}}%
+{h}(x_{{1}},x_{{2}},x_{4}),
\]
where ${h}(x_{{1}},x_{{2}},x_{4})$ is a smooth function. We substitute the
expression for $f$ into $\partial_{4}f=\tfrac{1}{2}{\mathrm{e}}^{-x_{2}}W^{2}$
and we get
\[
\partial_{4}{h}(x_{{1}},x_{{2}},x_{4})=\tfrac{1}{2}{\mathrm{e}}^{-x_{2}}%
W^{2}.
\]
By integrating the above equation, we obtain%
\[
{h}(x_{{1}},x_{{2}},x_{4})=\tfrac{1}{2}{\mathrm{e}}^{-x_{2}}W^{2}x_{4}%
+q(x_{1},x_{2}),
\]
for a smooth function $q$. Therefore,%
\[
f(x_{1},x_{2},x_{3},x_{4})=\tfrac{1}{2}{\mathrm{e}}^{-x_{1}}W^{1}x_{{3}%
}+\tfrac{1}{2}{\mathrm{e}}^{-x_{2}}W^{2}x_{4}+q(x_{1},x_{2}).
\]

We substitute the expressions for $\partial_{4}f$ and $\partial_{3}f$ into the
third and fourth equations of (\ref{GradientRS4}) and we get%
\[%
\begin{array}
[c]{ll}%
(2W^{1}-2W^{2}+\tfrac{4}{\mu})x_{3}-2{\mathrm{e}}^{x_{1}-x_{2}}W^{2}x_{4} &
\\[3pt]%
\multicolumn{1}{r}{+2{\mathrm{e}}^{x_{1}}\left(  -\partial_{1}q(x_{1}%
,x_{2})+2\partial_{2}q(x_{1},x_{2})-\mu\left(  W^{1}-\tfrac{1}{2}W^{2}\right)
\right) } & \!\!=0,\\[6pt]%
2{\mathrm{e}}^{x_{1}-x_{2}}W^{1}x_{3}+2\left(  W^{1}-W^{2}-\tfrac{2}{\mu
}\right)  x_{4}+3W^{3} & \\[3pt]%
\multicolumn{1}{r}{-2{\mathrm{e}}^{x_{2}}\left(  2\partial_{1}q(x_{1}%
,x_{2})-\partial_{2}q(x_{1},x_{2})-\mu\left(  \tfrac{1}{2}W^{1}-W^{2}\right)
\right) } & \!\!=0.
\end{array}
\]
Since the above equations must hold for all values of $x_{3}$ and $x_{4}$, we
easily get $W^{1}=W^{2}=0$. Hence, the above equations reduce to%
\[%
\begin{array}
[c]{l}%
{\tfrac{4}{\mu}x_{3}+2{\mathrm{e}}^{x_{1}}\left(  -\partial_{1}q(x_{1}%
,x_{2})+2\partial_{2}q(x_{1},x_{2})\right)  }=0,\\[6pt]%
{-\tfrac{4}{\mu}x_{4}+3W^{3}-2{\mathrm{e}}^{x_{2}}\left(  2\partial_{1}%
q(x_{1},x_{2})-\partial_{2}q(x_{1},x_{2})\right)  }=0.
\end{array}
\]
Clearly, the above equations cannot hold for all values of $x_{3}$ and $x_{4}%
$. Therefore, the solution of \eqref{RicciSoliton} we found for $g$ when
$\mu\neq0$ is not a gradient Ricci soliton.

\subsection{Case $\mu=0$}

We integrate equations (iii) and (v) of (\ref{System2}) with respect to
$B^{3}(x_{1},x_{2})$ and $B^{4}(x_{1},x_{2})$ respectively and we find
\begin{align}
B^{3}(x_{1},x_{2})  &  =e^{x_{1}}\left(  \tfrac{4}{3}-2e^{-x_{2}}\partial
_{1}B^{4}(x_{1},x_{2})\right)  ,\label{B3}\\[6pt]
B^{4}(x_{1},x_{2})  &  =F(x_{1})G(x_{2}), \label{B4}%
\end{align}
{ for some smooth functions} $F(x_{1})$ and $G(x_{2})$ satisfying%
\begin{align}
F^{\prime\prime}(x_{1})  &  =\tfrac{1}{3}cF(x_{1})-\tfrac{1}{3}F^{\prime
}(x_{1}),\label{EDOF}\\
G^{\prime}(x_{2})  &  =\tfrac{1}{2}cG(x_{2}), \label{EDOG}%
\end{align}
where $c$ is a real constant. By integrating (\ref{EDOF}) and (\ref{EDOG}) we
obtain:%
\begin{align}
F(x_{1})  &  =K_{1}e^{\left(  -\frac{1}{6}+\frac{1}{6}\sqrt{1+12c}\right)
x_{1}}+K_{2}e^{\left(  -\frac{1}{6}-\frac{1}{6}\sqrt{1+12c}\right)  x_{1}%
},\label{F}\\
G(x_{2})  &  =He^{\frac{1}{2}cx_{2}}, \label{G}%
\end{align}
for some constants $K_{i}$, $i=1,2$, and $H$. We substitute in (\ref{B4}) the
expressions (\ref{F}) and (\ref{G}) and we get
\begin{equation}
B^{4}(x_{1},x_{2})=\left(  K_{1}e^{\left(  -\frac{1}{6}+\frac{1}{6}%
\sqrt{1+12c}\right)  x_{1}}+K_{2}e^{\left(  -\frac{1}{6}-\frac{1}{6}%
\sqrt{1+12c}\right)  x_{1}}\right)  He^{\frac{1}{2}cx_{2}}. \label{B4b}%
\end{equation}
{Next, we substitute (\ref{B3}) and (\ref{B4b}) in equation (ix) of
(\ref{System2}) and we obtain}%
\begin{align*}
0  &  =H\left[  K_{1}\left(  \tfrac{3c}{4}+\left(  \tfrac{1}{6}-\tfrac{c}%
{4}\right)  \sqrt{1+12c}-\tfrac{1}{6}\right)  e^{\frac{x_{1}}{6}\sqrt{1+12c}%
}\right. \\
&  \left.  +K_{2}\left(  \tfrac{3c}{4}+\left(  \tfrac{c}{4}-\tfrac{1}%
{6}\right)  \sqrt{1+12c}-\tfrac{1}{6}\right)  e^{-\frac{x_{1}}{6}\sqrt{1+12c}%
}\right]  e^{\frac{cx_{2}-2x_{2}}{2}-\frac{x_{1}}{6}}-2,
\end{align*}
which clearly cannot hold for all values of $x_{1}$ and $x_{2}$. Therefore,
the case when $\mu=0$ does not correspond to a Ricci soliton.

The above results are summarized in the following.

\begin{theorem}
A four-dimensional proper generalized symmetric space of type $B$ (with the
invariant metric $g$ as described in \eqref{TypeB}) is a Ricci soliton if and
only if $\mu\neq0$. More precisely, for $\mu=0$, the metric $g$ described in
\eqref{TypeB} does not satisfy the Ricci soliton equation
\eqref{RicciSoliton}. For $\mu\neq0$, the Ricci soliton equation
\eqref{RicciSoliton} holds for the metric $g$ described in \eqref{TypeB}
together with a vector field
\begin{align*}
X  &  =\left(  -W^{1}+2W^{2}\right)  \left(  \frac{\partial}{\partial
x}+u\frac{\partial}{\partial u}\right)  +\left(  2W^{1}-W^{2}\right)  \left(
\frac{\partial}{\partial y}+v\frac{\partial}{\partial v}\right) \\
&  -\frac{4}{\mu}\left(  u\frac{\partial}{\partial u}+v\frac{\partial
}{\partial v}\right)  +W^{3}\frac{\partial}{\partial v},
\end{align*}
where $W^{1}$, $W^{2}$, $W^{3}$ are real constants. This Ricci soliton is not
a gradient one.
\end{theorem}

%\newpage

\end{document}